\documentclass[a4paper, 12pt]{article}
\usepackage{amssymb,amsbsy,float}
\usepackage{graphicx}
\usepackage{amsfonts}
\usepackage{latexsym}
\usepackage[centertags]{amsmath}
\usepackage{amsthm}
\usepackage{newlfont}
\newtheorem{theorem}{Theorem}

\newtheorem{claim}[theorem]{Claim}
\newtheorem*{claim*}{Claim}

\newtheorem{corollary}[theorem]{Corollary}

\newtheorem{lemma}[theorem]{Lemma}

\newtheorem{proposition}[theorem]{Proposition}

\newcommand{\finpreuve}{\ \rule{0.5em}{0.5em} \vspace{0.1 cm}}

\newcommand{\bx}{\mbox{$\textbf{x}$}}

\newcommand{\bA}{\mbox{$\textbf{A}$}}

\newcommand{\bb}{\mbox{$\textbf{b}$}}
\newcommand{\bz}{\mbox{$\textbf{z}$}}
\newcommand{\R}{\mbox{$\mathbb{R}$}}%
\newcommand{\N}{\mbox{$\mathbb{N}$}}%
\begin{document}
\title{Is having a unique equilibrium robust?}
\author{Yannick Viossat\thanks{E-mail:
yannick.viossat@polytechnique.org; address: CEREMADE, Universit\'e Paris-Dauphine, Place du Mar\'echal De Lattre De Tassigny,
75775 Paris Cedex 16, France; phone: 0033144054671; fax: 0033144054599.} \thanks{This paper was written
during my PhD, supervised by Sylvain Sorin, at the laboratoire d'\'econom\'etrie de l'Ecole polytechnique, Paris, and revised while at the Stockholm School of Economics. It was motivated by a
related note of Noa Nitzan (2005). It greatly benefited from
comments of Eilon Solan, Bernhard von Stengel, and participants to the game theory seminar of the Institut Henri Poincar\'e. Errors and shortcomings are mine.}} %
\date{\small{CEREMADE, Universit\'e Paris-Dauphine, Paris}} \maketitle

\begin{abstract}
We investigate whether having a unique equilibrium (or a given
number of equilibria) is robust to perturbation of the payoffs, both
for Nash equilibrium and correlated equilibrium. We show that the
set of $n$-player finite games with a unique correlated equilibrium
is open, while this is not true of Nash equilibrium for $n>2$. The
crucial lemma is that a unique correlated equilibrium is a
quasi-strict Nash equilibrium. Related results are studied. For
instance, we show that generic two-person zero-sum games have a
unique correlated equilibrium and that, while the set of symmetric
bimatrix games with a unique symmetric Nash equilibrium is not open,
the set of symmetric bimatrix games with a unique and quasi-strict symmetric Nash equilibrium is.\\

\noindent JEL classification C72\\

\end{abstract}

The relevance of a phenomenon arising in a game often hinges upon
this phenomenon being robust to perturbation of the payoffs of the
game. To establish such robustness results typically requires
proving that some of the properties of the game we initially
considered are themselves robust; that is, that the set of games
having these properties is open. We investigate here whether the set
of finite normal-form games with a unique equilibrium is open, for
Nash equilibrium, correlated equilibrium, and variants of Nash
equilibrium such as symmetric Nash equilibrium in symmetric games.

For Nash equilibrium in two-player games, the question has been solved by Jansen (1981),
who showed that the set of bimatrix games with a unique Nash
equilibrium is open. However, this result does not extend to
three-player games nor to symmetric equilibria of symmetric bimatrix
games: counterexamples are given in section \ref{UC-sec:remarks}.
Our main result is that, by contrast, for any number of players $n$,
the set of $n$-player finite games with a unique correlated
equilibrium is
open. 
This generalizes an earlier result of Noa Nitzan (2005).

An intuitive explanation of the discrepancy between the results on
Nash equilibrium and those on correlated equilibrium is as
follows: the proof of the openness of the set of bimatrix games
with a unique Nash equilibrium uses three ingredients :
upper-semi-continuity of the equilibrium correspondence, an
element of linearity in the structure of the set of equilibria
(the set of Nash equilibria of a bimatrix game is a finite union
of convex polytopes) and the fact that in bimatrix games, a unique
Nash equilibrium is quasi-strict (Jansen, 1981); that is, for each
player, no pure best response to the strategy of the other player
is outside the support of her own strategy.

For Nash equilibrium, the last two ingredients are specific to
two-player games, and the last one is also lacking for symmetric
Nash equilibria of symmetric bimatrix games. This accounts for our
negative results. By contrast, for any number of players, there is
strong element of linearity in the structure of the set of
correlated equilibria (this is a polytope). Furthermore, and this is
the crucial lemma, it may be shown that for any number of players, a
unique correlated equilibrium is a quasi-strict Nash equilibrium.
This follows from the strong complementary property of linear
programs and dual reduction arguments (Myerson, 1997).

The material is organized as follows: Definitions and notations are introduced in section $1$. Openness of
the set of games with a unique correlated equilibrium is proved in section 2. Section \ref{UC-sec:remarks}
groups remarks and related results: We first show that generic zero-sum games have a unique correlated
equilibrium and discuss the connections of our work with Nitzan's (2005). It is then shown that the set of
$3$-player games with a unique Nash equilibrium and the set of symmetric bimatrix games with a unique
symmetric Nash equilibrium are not open, but that the set of symmetric bimatrix games with a unique and
quasi-strict symmetric Nash equilibrium is. It is also shown that for any $n,k \geq 2$, the set of $n$-player
games with $k$ Nash equilibria (or $k$ extreme points to the set of correlated equilibria) is not open, though for
almost every game $G$, every game in a neighborhood of $G$ has the same number of equilibria as $G$. We then
study the structure of the set of Nash equilibria for bimatrix games at the relative boundary of the set of
games with a unique correlated equilibrium. Finally, we show that the set of two-player zero-sum games in
which one of the players has a unique optimal strategy is not open. This explains that the openness of the
set of games with a unique correlated equilibrium cannot be deduced easily from Hart and Schmeidler's (1989)
proof of existence of correlated equilibria.

\section{Definitions and main result}

Let $G$ be a finite $n$-player game. $I=\{1,2,...,n\}$ is the set
of players, $S^{i}$ the set of pure strategies of player $i$ and
$S^{-i}:=\times_{j \in I\backslash\{i\}} S^j$. The utility
function of player $i$ is $U^{i}:S=\times_{i \in I} S^{i} \to
\mathbb{R}$. As usual, $U^i$ is extended multilinearly to the set
of probability distributions over $S$. A pure strategy profile is
denoted by $s=(s^i, s^{-i})$ and a mixed strategy profile by
$\sigma=(\sigma^i,\sigma^{-i})$. The support of $\sigma^i$ is
denoted by%
$$Supp(\sigma^i):=\{s^i \in S^i: \sigma^i(s^i)>0\}$$
and the set of pure best-responses to $\sigma^{-i}$ by
$$PBR(\sigma^{-i}):=\{s^i \in S^i, \forall t^i \in S^i, U^i(s^i,\sigma^{-i}) \geq U^i(t^i,\sigma^{-i})\}$$
Finally, for any finite set $T$, the simplex of probability
distributions over $T$ is denoted by $\Delta(T)$.

A {\em correlated strategy} of the players in $I$ is a probability
distribution over the set $S$ of pure strategy profiles. Thus
$\mu=(\mu(s))_{s \in S}$ is a correlated strategy if:
\begin{equation}
\label{UC-eq:nnc}%
\mu(s)\geq 0 \hspace{5 mm} \forall s \in S
\end{equation}
\begin{equation}
\label{UC-eq:nmc}%
\sum_{s \in S} \mu(s)=1%
\end{equation}
Henceforth, the conditions in (\ref{UC-eq:nnc}) will be called \emph{nonnegativity constraints}. For $s^i$,
$t^i$ in $S^i$ and $\mu$ in $\Delta(S)$, let
$$h^{s^i,t^i}(\mu):=\sum_{s^{-i} \in S^{-i}} \mu(s)[U^{i}(s)-U^{i}(t^{i},s^{-i})]$$
where, as throughout, $s=(s^i, s^{-i})$.\\

\noindent \emph{Definition.} \textit{A correlated strategy $\mu$ is
a {\em correlated equilibrium} (Aumann, 1974) if it satisfies the
following {\em incentive constraints}:}
\begin{equation}
\label{UC-eq:incentives}%
h^{s^i, t^i}(\mu)
\geq 0, \hspace{5 mm} \forall i \in I, \forall s^i \in S^i, \forall t^{i} \in S^{i}%
\end{equation}%
Since conditions (\ref{UC-eq:nnc}), (\ref{UC-eq:nmc}) and
(\ref{UC-eq:incentives}) are all linear in $\mu$, it follows that
the set of correlated equilibria of any finite game is a polytope.

A $n$-player finite game has size $m_1 \times m_2 \times ...\times
m_n$ if, for every $i$ in $I$, the pure strategy set of player $i$
has cardinal $m_i$. Assimilating a game and its payoff matrices, a
$n$-player game of size $m_1 \times m_2\times ... \times m_n$ may be
seen as a point in $(\mathbb{R}^{m_1m_2...m_n})^n$, hence the
notions of a neighborhood of a game and of an open set of games. The
main result of this paper is that:
\begin{proposition}
\label{UC-prop:cor}%
The set of $n$-player games of size $m_1 \times m_2 \times ...
\times m_n$ with a unique correlated equilibrium is an open subset
of the set of games of size $m_1 \times m_2 \times
....\times m_n$. %
Furthermore, if a $n$-player finite game has a unique correlated
equilibrium $\sigma$, then the (unique) correlated equilibrium of
every nearby game has the same support as $\sigma$.
\end{proposition}
\section{Proof}
Let $G$ be a game with a unique correlated equilibrium and $(G_n)$
a sequence of games converging towards $G$. We need to show that,
for $n$ large enough, the game $G_n$ has a unique correlated
equilibrium. The proof runs as follows: Let $\sigma$ denote the
unique correlated equilibrium of $G$. A dual reduction argument
shows that $\sigma$ is a quasi-strict Nash equilibrium (lemma
\ref{lm:isquasistrict}). Together with the upper semi-continuity
of the Nash equilibrium correspondence this implies that, for $n$
large enough, $G_n$ has a quasi-strict Nash equilibrium with the
same support as $\sigma$ (lemma \ref{UC-lm:samesupport}). Since
two quasi-strict Nash equilibria with the same support satisfy the
same nonnegativity and incentive constraints with equality (lemma
\ref{lm:same-cons}), it follows that, for $n$ large enough, $G_n$
has a correlated equilibrium satisfying with equality the same
constraints as $\sigma$. Due to a general result on polytopes
(lemma \ref{Uniq-lm:polytope}), this implies that, for $n$
sufficiently large, the correlated equilibrium polytope of $G_n$
is a singleton. This completes the proof.

We begin with the result on polytopes: Let $(\bA_n)$ be a sequence
of $p \times q$ real matrices, $(\bb_n)$ a sequence of column
vectors of size $p$, and $J=\{1,...,p\}$. Let 
%
$$C_n=\{\bx \in \R^q, \bA_n\bx \geq \bb_n\}$$
($\bA_n\bx \geq \bb_n$ means that the weak inequality holds for
each coordinate). Make the following assumptions: first, $C_n$ is
uniformly bounded:
\begin{equation}
\label{Uniq-eq:unif-bounded} \exists M \in \R, \forall n \in \N,
\forall \bx \in C_n, \max_{j \in J} |(\bA_n\bx)_j| \leq M
\end{equation}
In particular, $C_n$ is a polytope. Second, $(\bA_n)$ and
$(\bb_n)$ converge respectively towards the matrix $\bA$ and the
vector $\bb$. Third, the ``limit polytope"
$$C=\{\bx \in \R^q, \bA\bx \geq \bb\}$$
is a singleton: $C=\{\bar{\bx}\}$. Let $J':=\{j \in J,
(\bA\bar{\bx})_j=b_j\}$ denote the set of constraints binding at
$\bar{\bx}$, and let $J"=J \backslash J'$. Finally, let $b_{n,j}$
denote the j$^{th}$ component of $\bb_n$.
\begin{lemma}
\label{Uniq-lm:polytope} If there exists $N \in \N$ such that, for
all $n \geq N$, there exists $\bx_n$ in $C_n$ with
\begin{equation}
\label{Uniq-eq:tightgivetight} \forall j \in J',
(\bA_n\bx_n)_j=b_{n,j}\end{equation}
then for $n$ large enough, $C_n$ is a singleton.
\end{lemma}
\begin{proof}
We begin by showing that, for $n$ large enough, all constraints
satisfied with strict inequality by $\bar{\bx}$ are satisfied with strict inequality by every $\bx$ in $C_n$: %
\begin{equation}
\label{eq:strictgivestrict} \exists N'\in \N, \forall n \geq N',
\forall \bx \in C_n, \forall j \in J", (\bA_n\bx)_j>b_{n,j}
\end{equation}
For all $n$ in $N$, let $\bz_n \in C_n$. To establish
(\ref{eq:strictgivestrict}), it suffices to show that for $n$
large enough, %
\begin{equation}
\label{eq:strictgive2} \forall j \in J", (\bA_n\bz_n)_i>b_{n,j}
\end{equation}
Due to (\ref{Uniq-eq:unif-bounded}), the sequence $(\bz_n)$ is
bounded. Furthermore, since $\bA_n \to \bA$ and $\bb_n \to \bb$, any
accumulation point $\bar{\bz}$ of $(\bz_n)$ satisfies $\bar{\bz} \in
C$, hence $\bar{\bz}=\bar{\bx}$. It follows that $(\bz_n)$ converges
to $\bar{\bx}$. Therefore $\bA_n\bz_n-\bb_n$ converges to $\bA\bar{\bx} -
\bb$. Since $(\bA\bar{\bx})_j - \bb_j>0$ for every $j$ in $J"$, it follows
that, for $n$ large enough, (\ref{eq:strictgive2}) is satisfied. This
completes the proof of (\ref{eq:strictgivestrict}).

We now prove the lemma: let $n \geq \max(N,N')$ and let $\bx_n$
satisfy (\ref{Uniq-eq:tightgivetight}). If $C_n$ is not a singleton,
then $C_n$ has an extreme point $\bz_n \neq \bx_n$. By basic
properties of extreme points, one of the constraints defining $C_n$
is binding at $\bz_n$ but not at $\bx_n$. Since $(\bA_n\bx_n)_j =b_{n,j}$
for all $j$ in $J'$, this implies that there exists $j$ in $J"$ such
that $(\bA_n\bz_n)_j=b_{n,j}$. This contradicts
(\ref{eq:strictgivestrict}). Therefore $C_n$ is a singleton.
\end{proof}

\begin{lemma}
\label{lm:isquasistrict} If an $n$-player finite game has a unique
correlated equilibrium $\sigma$ then this correlated equilibrium is
a quasi-strict Nash equilibrium.
\end{lemma}
\begin{proof}
The fact that $\sigma$ is a Nash equilibrium follows from the
existence of Nash equilibria and the fact that Nash equilibria are
correlated equilibria. What we really need to show is that $\sigma$ is quasi-strict. %

For each player $i$ in $I$, let $\alpha^i$ be a transition
probability over the set of pure strategies of player $i$: %
$$\begin{array}{cccc}
 \alpha^i: & S^i & \to & \Delta(S^i)\\
           & s^i & \to & \alpha^i \ast s^i
\end{array}
$$

A mixed strategy $\tau^i$ of player $i$ is
\emph{$\alpha^i$-invariant} if $\alpha^i \ast \tau^i =\tau^i$
where the mixed strategy $\alpha^i \ast \tau^i$ is defined by
$$\left[\alpha^i \ast \tau^i\right] (t^i)=\sum_{s^i \in S^i}
\tau^i(s^i) \, \left(\left[\alpha^i \ast s^i\right] (t^i)\right)
\hspace{1 cm} \forall t^i \in S^i$$
It follows from (Nau and McCardle, 1990, last paragraph of section 2
and proposition 2) that there exists a vector of transition
probabilities $\alpha$ such that, for every pure strategy profile
$s$ in $S$,
\begin{equation}
\label{eq:conddualvec} \sum_{i \in I} [U^i(\alpha^i \ast s^i,
s^{-i}) -U^i(s)] \geq 0
\end{equation}
with strict inequality if $s$ has probability zero in all correlated
equilibria. Fix such a vector $\alpha=(\alpha^i)_{i \in I}$. We
claim that:
\begin{claim}
\label{cl:invar} For every $i$ in $I$, the mixed strategy $\sigma^i$
is $\alpha^i$-invariant.
\end{claim}
\noindent This will be proved in the end. Assume that the pure
strategy $s^i$ does not belong to the support of $\sigma^i$ and let
$\tau=(s^i,\sigma^{-i})$. Since $\sigma$ is the unique correlated
equilibrium of $G$, it follows that every pure strategy profile $s$
in the support of $\tau$ has probability zero in all correlated
equilibria. Therefore, by definition of $\alpha$,
$$\tau(s)>0 \Rightarrow \sum_{k \in I} \left[U^k(\alpha^k \ast s^k,
s^{-k}) - U^k(s)\right] >0$$
It follows that
\begin{equation} \label{eq:arbitrage} \sum_{k \in I} \left[U^k(\alpha^k \ast \tau^k, \tau^{-k})-
U^k(\tau)\right] = \sum_{s \in S} \tau(s) \sum_{k \in I}
\left[U^k(\alpha^k \ast s^k, s^{-k}) - U^k(s)\right] >0
\end{equation}
(to prove the equality, use  $U^k(\alpha^k \ast \tau^k, \tau^{-k})= \sum_{s^k\in S^k} \tau^k(s^k)U^k(\alpha^k \ast s^k, \tau^{-k})$). 
For every $k \neq i$, $\tau^k=\sigma^k$ hence, by claim
\ref{cl:invar}, $\tau^k$ is $\alpha^k$-invariant. Therefore,
(\ref{eq:arbitrage}) boils down to $U^i(\alpha^i \ast \tau^i,
\tau^{-i})>U^i(\tau)$; that is,
$$U^i(\alpha^i \ast s^i, \sigma^{-i})>U^i(s^i, \sigma^{-i})$$
This implies that $s^i$ is not a best-response to $\sigma^{-i}$.
Since $s^i$ was an arbitrary strategy not in the support of
$\sigma^i$, it follows that $\sigma$ is quasi-strict.

It only remains to prove claim \ref{cl:invar}. The proof is based on
dual reduction (Myerson, 1997). Fix $\alpha$ as above. Note that,
due to (\ref{eq:conddualvec}), $\alpha$ is a dual vector in the
sense of Myerson (1997). Define the $\alpha$-reduced game $G/\alpha$
as in (Myerson, 1997). That is, in $G/\alpha$, the set of players
and the payoffs are as in $G$, but the mixed strategies available to
player $i$ are only those mixed strategies $\sigma^i$ of player $i$
in $G$ that are $\alpha^i$-invariant. Myerson (1997) shows that
$G/\alpha$ is a finite game, hence it has a Nash equilibrium. Let
$\tilde{\sigma}$ be a Nash equilibrium of $G/\alpha$. By definition
of $G/\alpha$, we may see $\tilde{\sigma}$ as a mixed strategy
profile of $G$, with $\tilde{\sigma}^i$ $\alpha^i$-invariant. Since
$\tilde{\sigma}$ is a Nash equilibrium of $G/\alpha$, it follows
from theorem 1 of Myerson (1997) that $\tilde{\sigma}$ is a Nash
equilibrium of $G$. Therefore $\tilde{\sigma}=\sigma$. Since
$\tilde{\sigma}^i$ is $\alpha^i$-invariant, this implies that
$\sigma^i$ is $\alpha^i$-invariant too. This completes the proof.
\end{proof}
The following lemma is a version of lemma 4.1 of Jansen (1981).
\begin{lemma}
\label{UC-lm:samesupport}%
If the $n$-player game $G$ has a unique Nash equilibrium $\sigma$
and this Nash equilibrium is quasi-strict, then there exists a
neighbourhood $\Omega_G$ of $G$ such that, for every game $\hat{G}$
in $\Omega_G$ and every Nash equilibrium $\hat{\sigma}$ of
$\hat{G}$, the support of $\hat{\sigma}$ is equal to the support of
$\sigma$ and $\hat{\sigma}$ is quasi-strict.
\end{lemma}
\begin{proof}
Let $(G_n)_{n \in \mathbb{N}}$ be a sequence of games converging
to $G$ and $\sigma_n$ a Nash equilibrium of $G_n$. To prove lemma
 \ref{UC-lm:samesupport}, it is enough to show that, for $n$ large enough,
 the support of $\sigma_n$ is equal to the support of $\sigma$ and $\sigma_n$ is
 quasi-strict. Since the Nash equilibrium correspondence is upper semi-continuous and
 since $G$ has a unique Nash equilibrium, it follows that $\sigma_n$ converges to
 $\sigma$. Therefore,
\begin{equation}
\label{eq:support} \exists N \in \N, \forall n \geq N, \forall i
\in I, Supp(\sigma^i) \subseteq Supp(\sigma_n^i)
\end{equation}
Furthermore, if $U^i(s^i,\sigma^{-i})<U^i(\sigma)$ then for $n$
large enough, $U_n^i(s^i,\sigma_n^{-i})<U_n^i(\sigma_n)$, where
$U_n^{i}$ denotes the utility function of player $i$ in the game
$G_n$. Therefore,
\begin{equation}
\label{eq:PBR} \exists N' \in \N, \forall n \geq N', \forall i \in
I, PBR(\sigma_n^{-i}) \subseteq PBR(\sigma^{-i})
\end{equation}
Finally, since the Nash equilibrium $\sigma$ is quasi-strict, it
follows that $PBR(\sigma^{-i})=Supp(\sigma^{-i})$. Together with
(\ref{eq:support}) and (\ref{eq:PBR}), this implies that for $n$
large enough:
$$Supp(\sigma_n^i) \subseteq PBR(\sigma_n^{-i}) \subseteq PBR(\sigma^{-i})=Supp(\sigma^i) \subseteq Supp(\sigma_n^{i})$$
Since the beginning and the end of this chain of inclusion are
equal, this is a chain of equality. In particular,
$Supp(\sigma_n^i) = PBR(\sigma_n^{-i}) = Supp(\sigma^i)$. The
result follows.\end{proof}
\begin{lemma}
\label{lm:same-cons} Let $G$ and $\hat{G}$ be two games with the
same set of players and strategies. Let $\sigma$ and
$\hat{\sigma}$ be Nash equilibria of, respectively, $G$ and
$\hat{G}$. Assume that $\sigma$ and $\hat{\sigma}$ have the same
support and are both quasi-strict. Then, among the nonnegativity
and incentive constraints defining correlated equilibria, $\sigma$
and $\hat{\sigma}$ satisfy the same constraints with
equality.\end{lemma}
\begin{proof}
Since, by assumption, $\sigma$ and $\hat{\sigma}$ have the same
support, the nonnegativity constraints they satisfy with equality
are the same. We now show that the incentive constraints they
satisfy with equality are also the same. Since $\sigma$ is a
product distribution, it follows that
$$h^{s^i,t^i}(\sigma)=\sigma^i(s^i)\left[U^i(s^i,\sigma^{-i})-U^i(t^i,\sigma^{-i})\right]
\hspace{1 cm} \forall i, \forall s^i, \forall t^i$$
Let $\tilde{S}=\times_i \tilde{S}^i$ denote the support of both
$\sigma$ and $\hat{\sigma}$. If $s^i \notin \tilde{S}^i$, then
$\sigma^i(s^i)=0$ hence $h^{s^i,t^i}(\sigma)=0$ for every $t^i$ in
$S^i$. If $s^i \in \tilde{S}^i$ and $t^i \in \tilde{S}^i$, then,
since $\sigma$ is a Nash equilibrium,
$U^i(s^i,\sigma^{-i})=U^i(t^i,\sigma^{-i})$ hence
$h^{s^i,t^i}(\sigma)=0$. Finally, if $s^i \in \tilde{S}^i$ and $t^i
\notin \tilde{S}^i$ then $\sigma^i(s^i)>0$ and, since $\sigma$ is
quasi-strict, $U^i(s^i,\sigma^{-i})-U^i(t^i,\sigma^{-i})>0$.
Therefore, $h^{s^i,t^i}(\sigma)>0$. Grouping these observations we
obtain that $h^{s^i,t^i}(\sigma)>0$ if and only if $s^i \in
\tilde{S}^i$ and $t^i \notin \tilde{S}^i$. The same result holds for
$\hat{\sigma}$ so that, letting $(\hat{h}^{s^i,t^i})_{s^i \in S^i,
t^i \in S^i}$ denote the linear forms associated with the correlated
equilibrium incentive constraints of $\hat{G}$, we have:
$$h^{s^i,t^i}(\sigma)=0 \Leftrightarrow \hat{h}^{s^i,t^i}(\hat{\sigma})=0$$
This completes the proof. %
\end{proof}

We now conclude. Let $G$ be a game with a unique correlated
equilibrium $\sigma$ and $(G_n)$ be a sequence of games converging
towards $G$. Let $C_n$ be the correlated equilibrium polytope of
$G_n$. Combining lemmas \ref{lm:isquasistrict},
\ref{UC-lm:samesupport} and \ref{lm:same-cons}, we obtain that, for
$n$ large enough, $G_n$ has a correlated equilibrium $\sigma_n$
satisfying with equality the same constraints as $\sigma$. By lemma
\ref{Uniq-lm:polytope}, this implies that for $n$ large enough,
$C_n$ is a singleton. This completes the proof of proposition
\ref{UC-prop:cor}.
\section{Remarks and related results}
\label{UC-sec:remarks}
\hspace{0.5 cm} 1. The fact that the set of $m_1 \times m_2 \times
... \times m_n$ games with a unique correlated equilibrium is
nonempty is obvious: any dominance solvable game has a unique
correlated equilibrium. Note also that generic two-player zero-sum
games have a unique correlated equilibrium. This follows from the
observation that: (i) a zero-sum game has a unique correlated
equilibrium if and only if it has a unique Nash equilibrium (Forges,
1990); (ii) generic zero-sum games have a unique Nash equilibrium
(Bohnenblust et al, 1950).

2. Nitzan (2005) proved independently and earlier a weaker version
of proposition \ref{UC-prop:cor}. More precisely, she proved that if
a two-player $m \times m$ game has a unique correlated equilibrium
and that this correlated equilibrium has full support, then every
nearby game has a unique correlated equilibrium and this correlated
equilibrium has full support. To prove this result with our method,
it suffices to note that if a game has a unique and completely mixed
Nash equilibrium, then every nearby game has a completely mixed Nash
equilibrium, and then to apply lemma \ref{Uniq-lm:polytope}. This
illustrates a difference between our arguments and Nitzan's: while
she uses a theorem of the alternative, we do not need any theorem of
the alternative to prove her results. (We do however use a theorem
of the alternative to prove proposition \ref{UC-prop:cor}. Indeed,
the proof of lemma \ref{lm:isquasistrict} uses Nau and McCardle's
(1990) characterization of strategy profiles with positive
probability in at least
one correlated equilibrium, which itself relies on a theorem of the alternative.)\\

3. The set of $3$-player games with a unique Nash equilibrium is
not open. The following $2 \times 2 \times 2$ counter-example is
adapted from (Flesch et al, 1997) and was provided by Eilon Solan
(personal communication).
\begin{equation}
\label{UC-eq:Eilongame}%
\left(\begin{array}{ll}
 1, 1, 1 & 0, 1, 1 \\
 1, 1, 0 & 1, 0, 1 \\
\end{array}\right)
\hspace{0.5 cm}
\left(\begin{array}{ll}
 1, 0, 1-\epsilon & 1, 1, 0 \\
 0, 1, 1          & 0, 0, 0 \\
\end{array}\right)
\end{equation}
%
Player $1$ chooses a row (Top or Bottom), player $2$ a column (Left
or Right) and player $3$ a matrix (West or East). For $\epsilon=0$,
there is a unique Nash equilibrium, in which all players play their
first strategy (this will be proved below). However, for
$\epsilon>0$, there is a continuum of Nash equilibria. Indeed, every
(partially) mixed strategy profile in which player $1$ plays Bottom
with probability less than $\epsilon/(1+\epsilon)$ and player $2$
and $3$ stick to their first strategy is a Nash equilibrium. Thus,
in $3$-player games, there are sequences of games with a continuum
of Nash equilibria converging towards a game with a unique Nash
equilibrium.

The game (\ref{UC-eq:Eilongame}) with $\epsilon=0$ also provides an example of a game with a unique Nash
equilibrium that is not quasi-strict. This calls for two remarks: first, while it is well known that
$3$-player games need not have a quasi-strict equilibrium, the counter-examples I found in the literature,
e.g. (Raghavan, 2002), are of games with several Nash equilibria. Thus, up to my knowledge, whether a unique
Nash equilibrium is necessarily quasi-strict was still open. Second, in two-player games, a unique Nash
equilibrium is necessarily quasi-strict, as shown by Jansen (1981), and as also follows from the fact that
every bimatrix game has a quasi-strict Nash equilibrium (Norde, 1999).\\

\noindent \emph{Proof that the game (\ref{UC-eq:Eilongame}) with
$\epsilon=0$ has a unique Nash equilibrium}: for $\epsilon=0$, the
game (\ref{UC-eq:Eilongame}) may be described as follows: player
$i+1$ (counted modulo $3$) wants to mismatch player $i$, except if
all players play their first strategy. Thus, in a hypothetical
equilibrium different from Top-Left-West, if $i$ plays in pure
strategy, then $i+1$ must mismatch $i$, $i+2$ mismatch $i+1$ and
$i+3=i$ mismatch $i+2$; therefore, $i$ must mismatch itself, a
contradiction. It follows that the only equilibrium in which one of
the players plays in pure strategy is Top-Left-West.

It remains to show that there are no completely mixed Nash
equilibria. By contradiction, let $x \in ]0,1[$ (resp. $y$, $z$) be
the probability of Bottom (resp. Right, East) in a hypothetical
completely mixed Nash equilibrium. Since player $1$ is indifferent
between Top and Bottom, we have $y(1-z) = (1-y)z + yz=z$, hence
$y>z$. Since the game is cyclically symmetric, it follows that
$y>z>x>y$, which cannot be. This completes the proof. \finpreuve \\

4. A corollary of lemma \ref{UC-lm:samesupport} is that:
\begin{corollary}\label{cor:unique-strict} The set of
games with a unique and strict Nash equilibrium is open.
\end{corollary}
Indeed, if a game has a unique and strict Nash equilibrium $\sigma$,
then by lemma \ref{UC-lm:samesupport}, every Nash equilibrium of
every nearby game has the same support as $\sigma$, hence is equal
to $\sigma$ as $\sigma$ is pure. I do not know whether the set of
games with a unique and
\emph{quasi}-strict Nash equilibrium is open.\\

5. The following example shows that, within the set of two-person
symmetric games, the set of games with a unique symmetric Nash
equilibrium is not open:
$$ \left(\begin{array}{ccc}
 -\epsilon, -\epsilon & 1, 0   & 1, 0   \\
 0, 1                 & 0, 0   & -1, -1 \\
 0, 1                 & -1, -1 & 0, 0   \\
\end{array}\right)
$$
For $\epsilon=0$, this game has a unique symmetric Nash
equilibrium: Top-Left. For $\epsilon>0$, it has $3$ symmetric Nash
equilibria: $(\frac{1}{1 + \epsilon}, \frac{\epsilon}{1 +
\epsilon}, 0)$, $(\frac{1}{1 + \epsilon}, 0, \frac{\epsilon}{1 +
\epsilon})$, $(\frac{3}{3 + 2\epsilon}, \frac{\epsilon}{3 +
2\epsilon}, \frac{\epsilon}{3 + 2\epsilon})$.
This is linked to the fact that, for $\epsilon=0$, the unique symmetric Nash equilibrium is not quasi-strict. Indeed, 
the openness of the set of bimatrix games with a unique Nash
equilibrium has the following analogue for symmetric games:
\begin{proposition} \label{prop:sym-quasistrict} Within the set of two-person symmetric games, the set of
two-person symmetric games with a unique symmetric Nash equilibrium and such that this Nash equilibrium is
quasi-strict is open.
\end{proposition}
\begin{proof}
Let $G$ be a two-person game with a unique symmetric Nash equilibrium, with support $\tilde{S}=\tilde{S^1}
\times \tilde{S^2}$ (To make things clear: there might be other, asymmetric Nash equilibria). Assume that the
unique symmetric equilibrium is quasi-strict. It follows from a variant of lemma \ref{UC-lm:samesupport}
that, within the set of two-person symmetric games, there exists a neighborhood $\Omega_G$ of $G$ such that,
for any game $G'$  in $\Omega_G$, any symmetric Nash equilibrium of $G'$ has support $\tilde{S}$ and is
quasi-strict.

Fix $G'$ in $\Omega_G$. Since $G'$ is symmetric, it has a symmetric
Nash equilibrium $\sigma$. To establish proposition
\ref{prop:sym-quasistrict}, it is enough to show that $G'$ has no
other symmetric Nash equilibrium. By contradiction, assume that $G'$
has a symmetric Nash equilibrium $\tau \neq \sigma$. For every
$\lambda$ in $\R$, define the symmetric mixed strategy profile
$\sigma_{\lambda}$ by $\sigma_{\lambda}^i=\lambda \tau^i +(1-
\lambda)\sigma^i$, for $i=1,2$.

There are five types of incentive and nonnegativity constraints that
$\sigma_{\lambda}$ must satisfy in order to be a (symmetric) Nash
equilibrium:

(i) $\sigma_{\lambda}(s) \geq 0, s \in \tilde{S}$;

(ii) $\sigma_{\lambda}(s) \geq 0, s \notin \tilde{S}$;

(iii) $h^{s^i,t^i}(\sigma_{\lambda}) \geq 0, s^i \in \tilde{S}^i,
t^i \notin \tilde{S}^i, i=1,2$;

(iv) $h^{s^i,t^i}(\sigma_{\lambda}) \geq 0, s^i \in \tilde{S}^i, t^i
\in \tilde{S}^i, i=1,2$;

(v)  $h^{s^i, t^i}(\sigma_{\lambda}) \geq 0, s^i \notin \tilde{S}^i,
t^i \in S^i, i=1,2$.

\noindent Using the fact that both $\sigma$ and $\tau$ are Nash
equilibria with support $\tilde{S}$, it is easily checked that for
every $\lambda$ in $\R$, $\sigma_{\lambda}$ satisfies (with
equality) all constraints of types (ii), (iv) and (v).

Moreover, since $\sigma \neq \tau$ and since the set of Nash
equilibria is compact, it follows that there exists a maximal value
of $\lambda$ such that $\sigma_{\lambda}$ is a (symmetric) Nash
equilibrium. Call this value $\lambda_{max}$. Since all symmetric
Nash equilibria of $G'$ have support $\tilde{S}$ and are
quasi-strict, they all satisfy with strict inequality all
constraints of types (i) and (iii), hence so does
$\sigma_{\lambda_{max}}$. Therefore, there exists
$\lambda>\lambda_{max}$ such that $\sigma_{\lambda}$ satisfies all
constraints of type (i) and (iii). Since, as mentioned in the
previous paragraph, $\sigma_{\lambda}$ also satisfies all other
constraints, it follows that $\sigma_{\lambda}$ is a (symmetric)
Nash equilibrium, contradicting the maximality of $\lambda_{max}$.
\end{proof}
Finally, a variant of the proof of corollary \ref{cor:unique-strict} shows that the set of $n$-player
symmetric games with a unique symmetric Nash equilibrium and such that this equilibrium is strict is open.\\

6. Until now, we focused on games with a unique equilibrium, but we
might also ask whether, for $k \geq 2$, having $k$ Nash equilibria,
or $k$ extreme points to the set of correlated equilibria, is a
robust property. The answer is negative:
\begin{proposition} For every $n, k \geq 2$, the set of $n$-player games
with $k$ ($\leq k$, $\geq k$) Nash equilibria is not open;
similarly, the set of $n$-player games with $k$ ($\leq k$, $\geq
k$) extreme points of the set of correlated equilibria is not
open.
\end{proposition}
We provide a counterexample for the case $n=2$, $k=3$. The
counter-example is easily generalized to any numbers $n,k \geq 2$.
For simplicity, call extreme correlated equilibria the extreme
points of the set of correlated equilibria. Consider the following
game:
\begin{equation}
\label{eq:EhudEilon}
\begin{array}{cc}
   & \begin{array}{ccc}
     L \hspace{0.7 cm} & M & \hspace{0.7 cm}R \\
    \end{array}     \\
  \begin{array}{c}
    T \\
    M \\
    B \\
  \end{array} & \left(\begin{array}{ccr}
                    0,0 & -1,-1 & -1,0 \\
                    -1,-1 & \epsilon,\epsilon & -1,0 \\
                    0,-1 & 0,-1 & 0,0 \\
                \end{array}\right) \\
 \end{array}
 \end{equation}
For $\epsilon=0$, both players can guarantee $0$ by playing the
third strategy, and this is the highest payoff they can get. It
follows that in any correlated equilibrium, the off-diagonal
strategy profiles have probability zero. Therefore, the three pure
strategy profiles on the diagonal are the only Nash equilibria of
the game and these are also the extreme correlated equilibria. For
$\epsilon<0$, the Nash equilibrium $(M,M)$ disappears, and only
two Nash equilibria remain, which again are also the extreme
correlated equilibria. For $\epsilon>0$, there are $4$ Nash
equilibria (resp. $6$ extreme correlated equilibria): $(T,L)$ and
the Nash equilibria (resp. extreme correlated equilibria) of the
$2 \times 2$ coordination game obtained by eliminating the
strategies $T$ and $L$.\\

7. The reason why a slight perturbation of the payoffs of
(\ref{eq:EhudEilon}) may alter the number of equilibria of the game
is that not all equilibria are quasi-strict. Indeed, Jansen (1981,
lemma 8.3 and remark 8.8) showed that if all equilibria of a
bimatrix game are quasi-strict, then (i) the game has a finite
number of equilibria and (ii) every nearby game has the same number
of equilibria.

Using the fact that a finite game has a finite number of equilibrium
components and upper-semi-continuity of the Nash equilibrium
correspondence, it is easy to generalize Jansen's result as follows:
if all equilibria of a $n$-player game are strongly stable (in the
sense of van Damme, 1991) then (i) and (ii) above hold. Since for
bimatrix games, all equilibria are strongly stable if and only if
all equilibria are quasi-strict, this indeed generalizes Jansen's
result. Furthermore, since for almost all games, all equilibria are
strongly stable (see, e.g., van Damme, 1991), it follows that (i)
and (ii) hold for almost all games.

Finally, it follows from the above discussion, and it is easy to
show directly that the set of $n$-player games with $k$ equilibria,
all strict, is open; but for $k \geq 2$, this set is actually void.
Indeed, if a game has $k \geq 2$ equilibria, all strict, then it
follows from an index argument that there exists at least $k-1$
mixed Nash equilibria, a contradiction. See corollary 2 and theorem
3 of Ritzberger (1994).\\

8. In bimatrix games, both the set of games with a unique Nash
equilibrium and the set of games with a unique correlated
equilibrium are open. Since there are games with a unique Nash
equilibrium but many correlated equilibria, the latter set is
included in the former. The following examples show that on the
relative boundary of the set of bimatrix games with a unique
correlated equilibrium, there are games with a continuum of Nash
equilibria, games with a finite number ($>1$) of Nash equilibria,
and games with a unique Nash equilibrium:
$$\left(\begin{array}{c}
    0     \\
 \epsilon \\
\end{array}\right)
\hspace{1 cm}
\left(\begin{array}{cc}
 1, 1 & 0, 0                 \\
 0, 0 & -\epsilon, - \epsilon \\
\end{array}\right)
\hspace{1 cm}
 \left(\begin{array}{lll|r}
 0, 0  & 2, 1  & 1, 2  & -1, x \\
 1, 2  & 0, 0  & 2, 1  & -1, x \\
 2, 1  & 1, 2  & 0, 0  & -1, x \\
 \hline
 x, -1 & x, -1 & x, -1 & 0, 0  \\
\end{array}\right)
$$
The left game is a one-person game. For $\epsilon>0$ it has a unique
correlated equilibrium. For $\epsilon=0$ it has a continuum of Nash
equilibria. The middle game has a unique correlated equilibrium
(Top-Left) for $\epsilon>0$, but two Nash equilibria for
$\epsilon=0$. The game on the right is adapted from (Nau and
McCardle, 1990, example 4). The $3 \times 3$ game in the top-left
corner is due to Moulin and Vial (1978). This $3 \times 3$ game has
a unique Nash equilibrium: $\left(1/3,1/3,1/3\right)$ for both
players, with payoff $1$; but putting probability $1/6$ on every
off-diagonal square yields a correlated equilibrium with payoff
$3/2$. Now consider the whole $4 \times 4$ game. For any value of
$x$, (4,4) is a Nash equilibrium. For $1<x \leq 3/2$, this is the
unique Nash equilibrium, but not the unique correlated equilibrium
(the correlated equilibrium with payoff 3/2 of the $3 \times 3$
top-left game induces a correlated equilibrium of the whole game).
For
$x>3/2$, this is the unique correlated equilibrium.\\

9. In order to prove the existence of correlated equilibria without
using a fixed point theorem, Hart and Schmeidler (1989) associate to
every finite game $G$ an auxiliary zero-sum game whose size depends
only on the size of $G$ and whose payoff matrix depends continuously
on the payoff matrices of $G$. In this auxiliary zero-sum game, the
optimal strategies of the maximizer correspond exactly to the
correlated equilibria of $G$, so that $G$ has a unique correlated
equilibrium if and only if, in the auxiliary game, the maximizer has
a unique optimal strategy. Therefore, to prove the openness of the
set of games with a unique correlated equilibrium, it would have
been enough to show that: \emph{If in a two-player zero-sum game,
one of the players has a unique optimal strategy, then in every
nearby zero-sum game this player has a unique optimal strategy.}

However, this turns out to be false: let $G_{\epsilon}$ be the
two-player zero-sum game with payoff matrix for the row player
$$
\begin{array}{cc}
                & \begin{array}{cc}
                         L & R \\
                  \end{array}                   \\
\begin{array}{c}
 T \\
 M \\
 B \\
\end{array}
                 & \left(\begin{array}{cc}
                         - \epsilon & 0  \\
                            0     & -1 \\
                            0     & -1 \\
                    \end{array}\right)          \\
\end{array}
$$
For $\epsilon=0$, the row player has a unique optimal strategy
(playing $T$). But for $\epsilon>0$, the row player has an infinite
number of optimal strategies: playing $T$ with probability
$1/(1+\epsilon)$ and playing $M$ and $B$ with any probabilities
summing to $\epsilon/(1+\epsilon)$.

\end{document}